\documentclass[english]{elsart}
\usepackage{babel}
\usepackage{epsfig,graphicx,times}
\usepackage{amssymb,amsmath}

\def\Ps{\mathcal{P}}

\newcommand{\Tsm}{\hspace*{0.6cm}}
\newcommand{\di}{m}

\def\@begintheorem#1#2{\list{}{\thm@body}%
  \item[]{\bf #1~#2.}\quad\it\ignorespaces}
\def\@opargbegintheorem#1#2#3{\list{}{\thm@body}%
  \item[]{\bf #1~#2~\ifrembrks #3\global\rembrksfalse\else (#3)\fi.}%
  \quad\it\ignorespaces}
\def\@endtheorem{\endlist}
\newtheorem{theorem}{Theorem}

\newtheorem{lemma}{Lemma}

\newtheorem{algo}{Algorithm}
\newtheorem{remark}{Remark}
\newtheorem{definition}{Definition}
\newtheorem{conjecture}{Conjecture}
\newcommand{\eop}{\hfill{$\Box$}}
\newcommand{\I}{\mathcal{I}}
\newcommand{\II}{{\!\!\overline{\I}(n,2)\!\!}}
\newcommand{\IC}{{\!\!\dot{\I}(n,2)\!\!}}
\newcommand{\Z}{\mathbb{Z}_n^\di}
\newcommand{\G}{\mathcal{G}}
\newenvironment{proof}
{\begin{trivlist}\item[]{{\sc Proof.}}}{\eop\noindent\end{trivlist}}

\begin{document}

\begin{frontmatter}
  \title{Integral point sets over $\mathbf{\Z}$}
  \author{Axel Kohnert}
  \ead{axel.kohnert@uni-bayreuth.de}
  \ead[url]{www.mathe2.uni-bayreuth.de}
  \address{University of Bayreuth, Department of Mathematics, D-95440 Bayreuth, Germany}
  \author{Sascha Kurz}
  \ead{sascha.kurz@uni-bayreuth.de}
  \ead[url]{www.wm.uni-bayreuth.de}
  \address{University of Bayreuth, Department of Mathematics, D-95440 Bayreuth, Germany}
  \begin{abstract}
    There are many papers studying properties of point sets in the Euclidean space $\mathbb{E}^\di$ or on integer
    grids $\mathbb{Z}^\di$, with pairwise integral or rational distances. In this article we consider the distances or
    coordinates of the point sets which instead of being integers are elements of $\mathbb{Z} /  \mathbb{Z}n$, and study the
    properties of the resulting combinatorial structures.
  \end{abstract}
\end{frontmatter}

\section{Introduction}

There are many papers studying properties of point sets in the Euclidean space $\mathbb{E}^\di$, with pairwise integral or rational distances (for short integral point sets or rational point sets, respectively), see \cite{integral_distances_in_point_sets} for an overview and applications. A recent collection of some classical open problems is given in \cite[Section 5.11]{research_problems}. Some authors also require that the points are located on an integer grid $\mathbb{Z}^\di$ \cite{Dimiev-Setting,cluster}. In this paper we modify the underlying space and study instead of $\mathbb{Z}$ the integers modulo $n$, which we denote by $\mathbb{Z}_n$. This was a suggestion of S. Dimiev. Our motivation was to gain some insight for the original problem in $\mathbb{Z}^\di$ and $\mathbb{E}^\di$. In the next subsection we shortly repeat the basic facts and questions about integral point sets in $\mathbb{Z}^\di$ and $\mathbb{E}^\di$.

\subsection{Integral point sets in $\mathbf{\mathbb{Z}^\di}$ and $\mathbf{\mathbb{E}^\di}$}
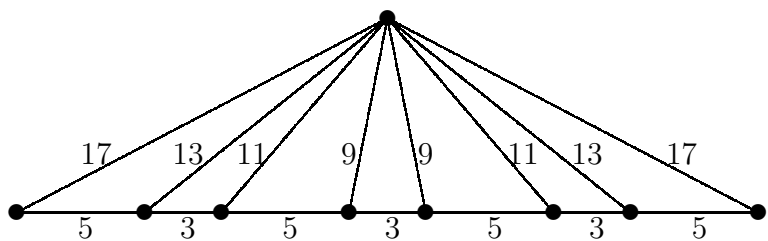
\begin{figure}[htp]
  \begin{center}
    \setlength{\unitlength}{0.34cm}
    \begin{picture}(29,8.8)
      \put(0,0.4){\line(1,0){29}}
      \put(0,0.4){\circle*{0.60}}
      \put(5,0.4){\circle*{0.60}}
      \put(8,0.4){\circle*{0.60}}
      \put(13,0.4){\circle*{0.60}}
      \put(16,0.4){\circle*{0.60}}
      \put(21,0.4){\circle*{0.60}}
      \put(24,0.4){\circle*{0.60}}
      \put(29,0.4){\circle*{0.60}}
      \put(14.5,8){\circle*{0.60}}
      \put(0,0.4){\qbezier(0,0)(7.25,3.8)(14.5,7.6)}
      \put(5,0.4){\qbezier(0,0)(4.75,3.8)(9.5,7.6)}
      \put(8,0.4){\qbezier(0,0)(3.25,3.8)(6.5,7.6)}
      \put(13,0.4){\qbezier(0,0)(0.75,3.8)(1.5,7.6)}
      \put(16,0.4){\qbezier(0,0)(-0.75,3.8)(-1.5,7.6)}
      \put(29,0.4){\qbezier(0,0)(-7.25,3.8)(-14.5,7.6)}
      \put(24,0.4){\qbezier(0,0)(-4.75,3.8)(-9.5,7.6)}
      \put(21,0.4){\qbezier(0,0)(-3.25,3.8)(-6.5,7.6)}
      \put(2.4,-0.6){$5$}
      \put(6.4,-0.6){$3$}
      \put(10.4,-0.6){$5$}
      \put(14.4,-0.6){$3$}
      \put(18.4,-0.6){$5$}
      \put(22.4,-0.6){$3$}
      \put(26.4,-0.6){$5$}
      \put(2.5,2.3){$17$}
      \put(6.1,2.3){$13$}
      \put(8.6,2.3){$11$}
      \put(12.7,2.3){$9$}
      \put(15.7,2.3){$9$}
      \put(19.2,2.3){$11$}
      \put(21.7,2.3){$13$}
      \put(25.4,2.3){$17$}
    \end{picture}
    \caption{A $2$-dimensional integral point set with $n=9$ and diameter $29$.}
    \label{fig_2_9_29}
  \end{center}
\end{figure}

So let us now consider integral point sets in $\mathbb{E}^\di$. If we denote the largest distance of an integral point set, consisting of $n$ points, as its diameter, the natural question for the minimum possible diameter $d(n,\di)$ arises, see Figure \ref{fig_2_9_29} for an example. Obviously we have $d(n,1)=n-1$. To avoid the corresponding trivial $1$-dimensional configuration in higher dimensions, it is common to request that an $\di$-dimensional integral point set is not contained in a hyperplane of $\mathbb{E}^\di$. We call a set of $\di+1$ points in $\mathbf{\mathbb{Z}^\di}$ or $\mathbf{\mathbb{E}^\di}$ degenerated, if the points are indeed contained in a hyperplane. There are quite a lot of constructions which show that $d(n,\di)$ exists for $n+1\ge \di$, see i.e. \cite{minimum_diameter}. Some exact values are determined in \cite{hab_kemnitz,phd_kurz,laue_kurz,sascha-alfred,dipl_piepmeyer}. The best known upper bound $d(n,\di)\in O\left(e^{c\,\log (n-\di) \log\log (n-!
 \di)}\right)$ is given in \cite{minimum_diameter}. 
For $\di=2$ Solymosi \cite{note_on_integral_distances} gives the best known lower bound $d(n,2)\ge cn$. For $\di=2$ and $n\ge 9$ the shape of the examples with minimum diameter is conjectured to consist of $n-1$ collinear 
points and one point apart \cite{sascha-alfred}, see Figure \ref{fig_2_9_29} for an example with $n=9$. We would like to remark that this conjecture is confirmed for $n\le 122$ by an exhaustive search \cite{sascha-alfred}. If for a fix $\rho>0$, we have a sequence of plane integral point set $\mathcal{P}_i$, each containing a collinear subset of cardinality least $n^\rho$, 
then the diameters of the $\mathcal{P}_i$ are in $\Omega\left(e^{c\,\log n\log\log n}\right)$ \cite{phd_kurz,sascha-alfred}. For $\di\ge 3$ we refer to \cite{phd_kurz,laue_kurz}, where some bounds and exact numbers are given.

Some authors require integral point sets to fulfill certain further conditions. The two classical conditions are, that no $\di+1$ points are contained in an $(\di-1)$-dimensional hyperplane, and that no $\di+2$ points are located on an $(\di-1)$-dimensional hypersphere. For ease of notation we speak of semi-general position in the first case and of general position if both conditions are fulfilled. We denote the minimum diameter of integral point sets in semi-general position by $\overline{d}(n,\di)$ and of integral point sets in general position by $\dot{d}(n,\di)$. For some small parameters the exact values have been determined in \cite{hab_kemnitz,kreisel,phd_kurz,sascha-alfred,dipl_piepmeyer}. 
We would like to remark that for dimension $m=2$ and $3\le n\le 36$ points, the examples with 
minimum possible diameter $\overline{d}(n,2)$, 
consist of points on a circle \cite{phd_kurz,sascha-alfred}.

A famous question of Erd\H{o}s asks for point sets in the plane with seven points in general position 
(i.e. no three on a line and nou four on a circle) with pairwise integral distances. 
Actually he first asked for such a set with five points, which was answered by Harborth \cite{0225.50001,0223.50008}, then for a set with six points, which was answered by Kemnitz \cite{hab_kemnitz}. Kemnitz even gives a construction for infinitely many such sets with coprime distances. For a long time no example consisting of seven points was known. Very recently one of the authors has discovered two such examples with diameters $22270$ and $66810$ \cite{kreisel}. For dimensions $m\ge 3$ we refer to \cite{phd_kurz,integral_distances_in_point_sets}.

As a specialization, integral point sets in general position, with all $n$ points on an integer grid $\mathbb{Z}^\di$, are called $n_\di$-clusters. Noll and Bell have found $n_\di$-clusters for $\di\le 5$ and $n\le\di+4$ but have no example for $n\ge\di+5$ \cite{cluster}. For $\di\ge 3$ even no integral point set in semi-general position with at least $\di+5$ points is known.

\begin{conjecture}{(Erd\H{o}s and Noll)}
  For any $\di>1$, $n>1$, there exists either none or an infinite number of non-isomorphic $n_\di$-clusters.
\end{conjecture}

An important invariant of an integral point set is its characteristic, which is defined as follows:

\begin{definition}
  Let $\mathcal{S}$ be a non-degenerated integral point set of $\di+1$ points in the $\di$-dimensional 
  Euclidean space $\mathbb{E}^\di$. By $V_m$ we denote the $\di$-dimensional volume of the simplex being
  formed by the convex hull of $\mathcal{S}$. Since the pairwise differences of $\mathcal{S}$ are integral
  and $\mathcal{S}$ is not degenerated we have $(V_m)^2\in\mathbb{N}\backslash\{0\}$. Thus $V_\di$ can be
  uniquely written as $V_\di=q\sqrt{c}$ with $q\in\mathbb{Q}$ and a squarefree integer $c$. This integer $c$
  is called the \textbf{characteristic} $\text{char}(\mathcal{S})$ of an integral simplex $\mathcal{S}$.
\end{definition}

The following theorem allows us to define the characteristic of an integral point set.

\begin{theorem}
  \label{thm_characteristic_old}
  In an $\di$-dimensional integral point set $\mathcal{P}$ each non-degenerate integral simplex $\mathcal{S}$
  has the same characteristic $\text{char}(\mathcal{S})$.
\end{theorem}

\begin{definition}
  Let $\mathcal{P}$ be an $\di$-dimensional integral point set and $\mathcal{S}\subseteq\mathcal{P}$ be an
  arbitrary $\di$-dimensional non-degenerate integral sub-simplex of $\mathcal{P}$. The \textbf{characteristic}
  $\text{char}(\mathcal{P})$ of $\mathcal{P}$ is given by $\text{char}(\mathcal{P})=\text{char}(\mathcal{S})$.
\end{definition}

For dimension $\di=2$ Theorem \ref{thm_characteristic_old} can be traced back at least to Kummer \cite{hab_kemnitz}, for $\di\ge 3$ we refer to \cite{paper_characteristic}. We would like to remark that if we are in the special case, where also the coordinates of an $\di$-dimensional integral point set $\mathcal{P}$ are integral, every subset $\mathcal{S}$ of $\mathcal{P}$, consisting of $\di+1$ points, has an integral volume. In our notation this means, that for an integral point set $\mathcal{P}$ in $\mathbb{Z}^m$ we have $\text{char}(\mathcal{P})=1$. So all $n_\di$-clusters have characteristic one. 

From  \cite{embedding,paper_characteristic} we know, that if $\mathcal{P}$ is an $\di$-dimensional integral point set in $\mathbb{E}^\di$ with characteristic $\text{char}(\mathcal{P})=1$, then there exists an embedding of $\mathcal{P}$ in $\mathbb{E}^\di$ using only rational coordinates. The existence of an embedding using only integral coordinates is an interesting open conjecture of \cite{embedding}.


\section{Integral point sets over $\mathbf{\Z}$}

In the previous section we have seen, that almost certainly there is a lot of hidden structure in the set of integral point sets which attain the minimum possible diameter and fulfill certain further conditions. Although the problem of integral point sets is a very classical one, not much progress has been achieved towards structure results or tight bounds on the minimum diameter. The idea of this paper is to study similar problems, which might be easier to handle, but may give some insight in the original problem. At first we want to consider the study of integral point sets in $\mathbb{Z}^\di$ as our \textit{original problem} and relate it to some problem of point sets in $\mathbb{Z}_n^\di$.

So let $\mathcal{P}'$ be an integral point set over $\mathbb{Z}^\di$. To relate $\mathcal{P}'$ 
to a set $\mathcal{P}$ of points in $\mathbb{Z}_n^\di$ we consider 
the canonical  mapping $\phi_n:\mathbb{Z}\rightarrow\mathbb{Z}_n$, $x\mapsto x+\mathbb{Z}n=\overline{x}$, which maps coordinates in $\mathbb{Z}^m$ to coordinates in $\Z$. If $n$ is suitably
large no two points of $\mathcal{P}'$ will be mapped onto the same point in $\mathcal{P}$. To be able to translate results in $\mathbb{Z}_n^\di$ back to $\mathbb{Z}^\di$, we define the inverse mapping $\Psi_n:\mathbb{Z}_n\rightarrow\{0,\dots,n-1\}$ by $\Psi(\phi_n(x))=x$ for $x\in\{0,\dots,n-1\}$. As an abbreviation we set $\Psi_n(x)=\widehat{x}$ and $\phi_n(x)=\overline{x}$, whenever the value of $n$ is clear from the context. Since points in $\mathcal{P}'$ have integral distances in $\mathbb{Z}^\di$ we need a similar definition of integral distances in $\mathbb{Z}_n^\di$. The most natural way to define an integral distance over $\Z$ is:

\begin{definition}
  Two points $(u_1,\dots,u_\di),(v_1,\dots,v_\di)\in \Z$ are at \textbf{integral distance}, if
  there exists a number $d\in\mathbb{Z}_n$ with
  $$
    \sum_{i=1}^{\di}(u_i-v_i)^2=d^2.
  $$
\end{definition}

\noindent
With this definition an integral point set $\mathcal{P}'$ over $\mathbb{Z}^\di$ is mapped via $\phi_n$ onto an integral point set $\mathcal{P}$ over $\Z$. Since $\phi_n$ may map some point set $\mathcal{P}'$ over $\mathbb{Z}^\di$, which is not contained in a hyperplane of $\mathbb{Z}^\di$, onto a point set $\mathcal{P}'$, where all points are contained in a hyperplane of $\Z$, we do not make any requirements on the distribution of the points in an integral point set over $\Z$ in the first run. The next definition to translate from $\mathbb{Z}^\di$ or $\mathbb{E}^\di$ to $\mathbb{Z}_n^\di$ is the minimum diameter. In $\mathbb{Z}^\di$ and $\mathbb{E}^\di$ we need the concept of a minimum diameter to get a finite space, whereas $\Z$ is finite for itself. So we find it natural to consider the maximum number of integral points.

\begin{definition}
  By $\I(n,\di)$ we denote the maximum number of points in $\Z$ with pairwise integral distances.
\end{definition}

\begin{theorem}
  \label{thm_trivial}
  $\I(n,1)=n$, $\I(1,\di)=1$, and $\I(2,\di)=2^\di$.
\end{theorem}
\begin{proof}
  Because there are only $n^\di$ different elements in $\Z$ we have the trivial upper bound $\I(n,\di)\le n^\di$. 
  This upper bound is only attained if $\di=1$ or $n\le 2$, since $\mathbb{Z}_n$ has at least one quadratic non
  residue for $n\ge 3$.
\end{proof}

\begin{table}[ht]
  \scriptsize
  \begin{center}
    \begin{tabular}{|r|r|r|r|r|r|r|r|r|r|r|r|r|r|}
      \hline
      \!\!\!$\di \backslash n$\!\!\!& 3 & 4 & 5 & 7 & 8 & 9 & 11 & 13 & 16 & 17\\ 
      \hline
      2&3&8&5&7&16&27&11&13&64&17\\ 
      3&4&16&25&8&64&81&11&169&256&289\\ 
      4&9&32&25&49&512&324&121&$\ge 169$&1024&\\ 
      5&27&128&125&343&2048&$\ge 893$&$\ge 1331$&$\ge 2197$&&\\ 
      6&33&256&$\ge 125$&&$\ge 15296$&&&&&\\ 
      7&$\ge 35$&1024&&&$\ge 81792$&&&&&\\   
      \hline
    \end{tabular} 
    \caption{Values of $\mathcal{I}(n,\di)$ for small parameters $n$ and $\di$.}
    \label{table_val_n_m}
  \end{center}
\end{table}

\noindent
For $n\ge 3$ we so far were not able to derive explicit formulas for $\I(n,\di)$ and so we give in Table \ref{table_val_n_m} some values for small parameters $n$ and $\di$, obtained by exhaustive enumeration via clique search, which we will describe in the next subsection. Further exact values or lower bounds can be determined using Theorem \ref{thm_trivial} and \ref{thm_cartesian} of Subsection \ref{sub_sec_hom}.

\subsection{Exhaustive enumeration of integral point sets over $\mathbf{\Z}$ via clique search}

In this subsection we describe how the exact values $\mathcal{I}(n,\di)$ of Table \ref{table_val_n_m} 
were obtained. 
We model our problem as a graph $\G$, so that the cliques (i.e. complete subgraphs) of $\G$ are in bijection to integral point sets over $\Z$. Therefore we choose the elements of $\Z$ as vertices and connect $x,y\in\Z$ via an edge, if and only if $x$ and $y$ are at integral distance.

To determine $\I(n,\di)$, we only have to determine the maximum cardinality of a clique of $\G$. Unfortunately this is an $\mathcal{NP}$-hard problem in general, but practically this approach was also successful in the case of integral point sets over $\mathbb{E}^\di$ \cite{phd_kurz,sascha-alfred}, due to good heuristic maximum-clique algorithms. Besides an implementation of the Bron-Kerbosch algorithm \cite{0261.68018} 
written by ourself we use the software 
package \textsc{cliquer} \cite{cliquer,1019.05054} of Niskanen and \"Osterg{\aa}rd.  

By prescribing points or distances of an integral point set $\Ps$, it is possible to reduce the complexity for the clique-search algorithm. The first variant is, that due to symmetry we can assume that the point $0=(\overline{0},\dots,\overline{0})\in\Z$ is part of $\Ps$. As vertices of $\G$ we choose the points in $\Z\backslash\{0\}$, which have an integral distance to $0$. Again two vertices $x,y\in\G$ are joined by an edge, if the corresponding points are at integral distance.

For the second variant we consider the set $D_{n,m}$ of all points $d=(d_1,\dots,d_\di)\in\Z$, which
have an integral distance to $0$ and which fulfill $\widehat{d}_i\le\left\lfloor\frac{n}{2}\right\rfloor$, for all $1\le i\le m$. So for every two points $u=(u_1,\dots,u_\di)\neq v=(v_1,\dots,v_d)\in\Z$, having an integral distance,
the tuple
$$
  \delta_n(u,v)=\left(\,\overline{\min\left(|\widehat{u}_1-\widehat{v}_1|,n-|\widehat{u}_1-\widehat{v}_1|\right)},
  \dots,\overline{\min(|\widehat{u}_\di-\widehat{v}_\di|,n-|\widehat{u}_\di-\widehat{v}_\di|)}\,\right)
$$ 
is an element of $D_{n,m}$. Actually we consider the vector of the Lee weights \cite{0907.94001} 
of the coordinates of the difference $u-v$. 
Now we choose an arbitrary numbering of this set $D_{n,m}\backslash\{0\}=\{{e_0},\dots,{e_{|D_{n,m}|-2}}\}$ and consider 
the graphs $\G_i$, which consist of the points of $\Z\backslash\{0,{e_i}\}$, with integral distances to $0$ and $e_i$, as vertices. Two vertices $x\neq y\in\G$ are joined by an edge if and only if the corresponding points fullfill $\delta_n(x,y)={e_j}$ with $i\le j$. Again one can show, that an integral point set in $\Z$ corresponds to a clique in
some graph $\G_i$ and vice versa. For some values of $n$ and $\di$ it is worth to put some effort 
in a suitable choice of the numbering of $D_{n,m}\backslash\{0\}$.

\subsection{Hamming spaces and homomorphisms}
\label{sub_sec_hom}
In this subsection we want to relate the problem of integral point sets over $\Z$ to problems in Hamming spaces.
In coding theory the Hamming distance $h(u,v)$ of two vectors $u=(u_1,\dots,u_\di)$, $v=(v_1,\dots,v_\di)$ $\in\Z$ is the number of positions $i$ where $u_i$ and $v_i$ differ. Normally one is interested in large subsets of $\Z$ where all the Hamming distances are either $0$ or larger than a given constant $c$. In our subject, we are interested in large subsets of $\Z$, where all the Hamming distances are taken from a specific proper subset of $\{0,1,\dots,\di\}$. This point of view has been proven useful i.e. also in the $0/1$-Borsuk problem in low dimensions, see \cite{ziegler}. Here we also want to mention the study of two-weight codes, see i.e. \cite{calderbank-kantor-two-weight-86,1268198}.

So let us go back to the determinantion of $\I(n,\di)$. As there are trivial formulas for $\I(1,\di)$ 
and $\I(2,\di)$, the next open case for fixed ring order 
$n$ is the determination of $\I(3,\di)$. Due to $1^2\equiv 2^2\equiv 1\mod 3$, integral point sets over $\mathbb{Z}_3^\di$ correspond to sets of $\mathbb{Z}_3^\di$ with Hamming distances $h(u,v)\not\equiv 2\mod 3$. So this is our first example of a selection problem in a Hamming space.

For the determination of $\I(2n,\di)$ we can utilize homomorphisms to make the problem easier. Therefore we need some definitions.

\begin{definition}
  For an integer $n$ we define the mapping $\tilde{\varphi}_{2n}:\mathbb{Z}_{2n}\rightarrow\mathbb{Z}_n$,
  $x\mapsto \widehat{x} +\mathbb{Z}n$, and by $\varphi_{2n,\di}$ we denote its extensions to $\mathbb{Z}_{2n}^\di$.
\end{definition}

\begin{definition}
  The weight function $\tilde{w}_{2n}:\mathbb{Z}_n^2\rightarrow\mathbb{Z}_{2n}$ is defined by
  $(u_i,v_i)\mapsto (\widehat{u}_i-\widehat{v}_i)^2+\mathbb{Z}\cdot 2n$.
  $$
    \mathbb{H}_{2n}^m:=\left\{S\subseteq\mathbb{Z}_n^m\mid \forall\,s_1,s_2\in S:\,\exists d\in\mathbb{Z}_{2n}:\, 
    d^2=w(s_1,s_2)\right\},
  $$
  where $w_{2n,\di}:\left(\mathbb{Z}_n^\di\right)^2\rightarrow\mathbb{Z}_{2n}$ is given by 
  $\left((u_1,\dots,u_\di),(v_1,\dots,v_\di)\right)\mapsto\sum\limits_{i=1}^\di\tilde{w}_{2n}(u_i,v_i)$. By 
  $\mathbb{I}_n^\di$ we denote the set of integral point sets in $\mathbb{Z}_n^\di$.
\end{definition}

\begin{lemma}
  \label{lemma_divisibility}
  $$2^\di\,|\,\I(2n,\di).$$  
\end{lemma}
\begin{proof}
  We consider the ring homomorphism $\varphi_{2n,\di}$ and restrict it to 
  $\varphi'_{2n,\di}:\mathbb{I}_{2n}^{\di}\rightarrow\mathbb{H}_{2n}^\di$. If $\Ps$ is an element
  of $\mathbb{H}_{2n}^\di$ then the preimage $\varphi_{2n,\di}^{-1}(\Ps)$ is an integral point set, due to
  $(x+n)^2\equiv x^2+n\mod 2n$ for odd $n$ and $(x+n)^2\equiv x^2\mod 2n$ for even $n$. For all
  $x\in\mathbb{Z}_n^\di$ we have $|\varphi_{2n,\di}^{-1}(x)|=2^\di$.
\end{proof}

So for the determination of $\I(2n,\di)$, it suffices to determine the maximum cardinality of the elements of $\mathbb{H}_{2n}^\di$, which actually are subsets of $\mathbb{Z}_n^\di$.
$$
 \I(2n,\di)=2^m\cdot \max_{S\in\mathbb{H}_{2n}^\di} |S|
$$
As an example we want to apply this result for $n=2$. Here $w_{4,\di}$ is exactly the Hamming distance in $\mathbb{Z}_2^\di$. Since the squares of $\mathbb{Z}_4$ are given by $\{0,1\}$, we conclude that $\mathbb{H}_4^\di$ is the set of all subsets of $\mathbb{Z}_2^\di$, with Hamming distance congruent to $0$ or $1$ modulo $4$. With the mapping $\varphi'_{4,\di}$ at hand, we can exhaustively generate the maximal sets in $\mathbb{H}_{4}^\di$, via a clique search, to extend Table \ref{table_val_n_m}:
$$
  \left(\I(4,\di)\right)_{\di\le 12}=4,8,16,32,128,256,1024,4096,16384,32768,65536,131072. 
$$

The next theorem shows, that it suffices to determine $\I(a,\di)$ for prime powers $a=p^r$.
\begin{theorem}
  \label{thm_cartesian}
  For two coprime integers $a$ and $b$ we have $\I(a\cdot b,\di)=\I(a,\di)\cdot\I(b,\di)$.
\end{theorem}
\begin{proof}
  Since $a$ and $b$ are coprime we have $\mathbb{Z}_{ab}\simeq\mathbb{Z}_a\times\mathbb{Z}_b$. If $\Ps$ is an
  integral point set in $\mathbb{Z}_a\times\mathbb{Z}_b$, then the projections into $\mathbb{Z}_a$ and $\mathbb{Z}_b$
  are also integral point sets. If on the other hand, $\mathcal{P}_1$ and $\mathcal{P}_2$ are integral point sets 
  over $\mathbb{Z}_a$ and $\mathbb{Z}_b$, respectively, then $\Ps:=\mathcal{P}_1\times\mathcal{P}_2$ is an integral
  point set over $\mathbb{Z}_a\times\mathbb{Z}_b$, due to a straight forward calculation.
\end{proof}

If we drop the condition that $a$ and $b$ are coprime 
Theorem \ref{thm_cartesian} does not remain valid in general. 
One can see this by looking  at the example $\I(2,3)\cdot\I(4,3)>\I(8,3)$ in table \ref{table_val_n_m}. 
Also $\I(a,\di)\mid \I(a\cdot b,\di)$ does not hold in general, as on can see by a look at the example $\I(3,3)\!\nmid\,\I(9,3)$. We would like to mention, that in a recent preprint \cite{finite_fields} the exact values of $\I(p,2)$ and $\I(p^2,2)$ have been determined.

\begin{theorem}
  For a prime $p\ge 3$ we have
  $$
    \I(p,2)=p\quad\text{and}\quad\I(p^2,2)=p^3.
  $$
\end{theorem}

\subsection{Integral point sets over the plane $\mathbf{\mathbb{Z}_n^2}$}

In Theorem \ref{thm_trivial} we have given an exact formula for $\I(n,1)$. So, if we fix the dimension $\di$, the next case is the determination of $\I(n,2)$. At first we give two constructions to obtain lower bounds for $\I(n,2)$.


\begin{lemma}
  \label{construction_1}
  If the prime factorization of $n$ is given by $n=\prod\limits_{i=1}^s p_i^{r_i}$, with pairwise different primes
  $p_i$, we have
  $$
    \I(n,2)\ge n\cdot\prod_{i=1}^s p_i^{\left\lfloor\frac{r_i}{2}\right\rfloor}.
  $$
\end{lemma}
\begin{proof}
  We choose the points $(u_i,v_j\overline{k})$, where $u_i,v_j\in\mathbb{Z}_n$ and
  $k=\prod\limits_{i=1}^s p_i^{\left\lceil\frac{r_i}{2}\right\rceil}$.
  Since
  $$
    (u_{i_1}-u_{i_2})^2+(v_{j_1}\overline{k}-v_{j_2}\overline{k})^2=(u_{i_1}-u_{i_2})^2,
  $$
  all occurring distances are integral.
\end{proof}

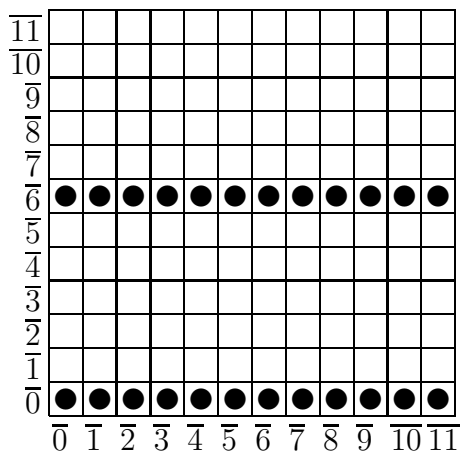
\begin{figure}[htp]
  \begin{center}
    \setlength{\unitlength}{0.45cm}
    \begin{picture}(13.2,13.2)
      \multiput(1.15,1)(0,1){13}{\line(1,0){12}}
      \multiput(1.15,1)(1,0){13}{\line(0,1){12}}
      \multiput(1.15,1)(1,0){12}{\put(0.5,0.5){\circle*{0.6}}}
      \multiput(1.15,1)(1,0){12}{\put(0.5,6.5){\circle*{0.6}}}
      \put(1.25,0){$\overline{0}$}
      \put(2.25,0){$\overline{1}$}
      \put(3.25,0){$\overline{2}$}
      \put(4.25,0){$\overline{3}$}
      \put(5.25,0){$\overline{4}$}
      \put(6.25,0){$\overline{5}$}
      \put(7.25,0){$\overline{6}$}
      \put(8.25,0){$\overline{7}$}
      \put(9.25,0){$\overline{8}$}
      \put(10.25,0){$\overline{9}$}
      \put(11.25,0){$\overline{10}$}
      \put(12.35,0){$\overline{11}$}
      \put(0.45,1.1){$\overline{0}$}
      \put(0.45,2.1){$\overline{1}$}
      \put(0.45,3.1){$\overline{2}$}
      \put(0.45,4.1){$\overline{3}$}
      \put(0.45,5.1){$\overline{4}$}
      \put(0.45,6.1){$\overline{5}$}
      \put(0.45,7.1){$\overline{6}$}
      \put(0.45,8.1){$\overline{7}$}
      \put(0.45,9.1){$\overline{8}$}
      \put(0.45,10.1){$\overline{9}$}
      \put(0,11.1){$\overline{10}$}
      \put(0,12.1){$\overline{11}$}
    \end{picture}
  \end{center}
  \caption{An integral pointset over $\mathbb{Z}_{12}^2$ constructed via Lemma \ref{construction_1}.}
  \label{fig_construction_1}
\end{figure}

An example of the construction of Lemma \ref{construction_1} is given in Figure \ref{fig_construction_1}, for $n=12=2^2\cdot 3$.

In the case of $n=2 \mod 4$ we can improve the above lemma:

\begin{lemma}
  \label{construction_2}
  If the prime factorization of $n$ is given by $n=2\cdot\prod\limits_{i=2}^s p_i^{r_i}$, with
  pairwise different primes $p_i\neq 2$ we have
  $$\I(n,2)\ge 2n\cdot\prod_{i=2}^s p_i^{\left\lfloor\frac{r_i}{2}\right\rfloor}.$$
\end{lemma}
\begin{proof}
  We choose the points $(u_i,v_j\overline{k})$, where $u_i,v_j\in\mathbb{Z}_n$ and
  $k=\prod\limits_{i=2}^s p_i^{\left\lceil\frac{r_i}{2}\right\rceil}$. Since $2k^2\equiv 0\mod n$ and
  $$
    (u_{i_1}-u_{i_2})^2+(v_{j_1}\overline{k}-v_{j_2}\overline{k})^2= (u_{i_1}-u_{i_2})^2
    + (v_{j_1}^2+v_{j_2}^2)\overline{k^2}
  $$
  either
  $$
    (u_{i_1}-u_{i_2})^2+(v_{j_1}\overline{k}-v_{j_2}\overline{k})^2= (u_{i_1}-u_{i_2})^2
  $$
  or
  $$
    (u_{i_1}-u_{i_2})^2+(v_{j_1}\overline{k}-v_{j_2}\overline{k})^2=(u_{i_1}-u_{i_2}+\overline{k^2})^2
  $$
  holds.
\end{proof}

\begin{conjecture}
  \label{conj_kohnert_kurz}
  For all $n\in\mathbb{N}$ either the lower of Lemma \ref{construction_1} or the lower bound of
  Lemma \ref{construction_2} is tight.
\end{conjecture}

\begin{remark}
  By Theorem \ref{thm_cartesian} and an exhaustive enumeration of integral point sets over $\mathbb{Z}_n^2$, via
  clique search, we have verified Conjecture \ref{conj_kohnert_kurz} up to $n=307$.
\end{remark}

If $n$ is squarefree and $2$ does not divide $n$, then our constructions from Lemma \ref{construction_1} and
Lemma \ref{construction_2} yield point sets of the form $\mathcal{P}=\{(u,0)\mid u\in\mathbb{Z}_n\}$. This is somewhat similar to the situation in $\mathbb{E}^2$, where integral collinear point sets with small diameter can consist of many points. Since we also want to speak of collinear point sets in $\mathbb{Z}_n^2$ we give:

\begin{definition}
  \label{definition_collinear}
  A set of $r$ points $(u_i,v_i)\in\mathbb{Z}_n^2$ is collinear, if there are $a,b,t_1,t_2,w_i\in \mathbb{Z}_n$ 
  with
  $$
    a+w_it_1=u_i\,\,\quad\text{and}\,\,\quad b+w_it_2=v_i.
  $$
\end{definition}

Let us first look at collinearity from the algorithmic point of view. Checking three points for being collinear, by running through the possible values of $a,b,t_1,t_2,w_i\in \mathbb{Z}_n$, would cost $\mathcal{O}(n^7)$ time. Setting, w.l.o.g., $a=u_1$, $b=v_1$, $w_1=\overline{0}$ reduces this to $\mathcal{O}(n^4)$.
If $n$ is prime, then we are working in a field, and there is an easy and well known way to check, whether three points are collinear, in $\mathcal{O}(1)$ time:

\begin{lemma}
  \label{lemma_characterization_collinear}
  For a prime $n$ the points $(u_1,v_1),(u_2,v_2),(u_3,v_3)\in\mathbb{Z}_n^2$ are collinear, if and only if
  \begin{eqnarray}
    \label{eqn_collinear}
    \left|\begin{array}{ccc}u_1&v_1&\overline{1}\\u_2&v_2&\overline{1}\\u_3&v_3&\overline{1}\\\end{array}\right|=\overline{0}.
  \end{eqnarray}
\end{lemma}


We remark that in $\mathbb{Z}_8$ the points $(\overline{0},\overline{0})$, $(\overline{2},\overline{4})$, 
$(\overline{4},\overline{4})$ fulfill equation (\ref{eqn_collinear}), but are not collinear with respect to Definition
\ref{definition_collinear}. So in general equation (\ref{eqn_collinear}) is necessary but not sufficient for three points to be collinear. We propose the development of a fast algorithm, which checks three points in $\mathbb{Z}_n^2$ for being collinear, as an interesting open problem. In practice one simply determines for each pair $x,y\in\mathbb{Z}_n^2$, whether the triple $0,x,y$ is collinear or not, in a precalculation.

The study of collinear point sets is motivated by the situation in the case of non-modular point sets.
Due to a theorem of Erd\H{o}s each integral point set  in  $\mathbb{E}^2$, with infinitely many points,
is located on a line \cite{ErdoesAnning1,ErdoesAnning2}. And, as already mentioned in the introduction the, non-collinear integral point sets in $\mathbb{E}^2$ with minimum diameter, are conjectured to consist of $n-1$ collinear points and one point apart.

In this context we would like to mention a theorem, which was recently proven in \cite{finite_fields}.
\begin{theorem}
  \label{conj_collinear}
  For $p$ being a prime, with $p\equiv 3\mod 4$, each integral point set over $\mathbb{Z}_p^2$, consisting of $p$
  points, is collinear.
\end{theorem}

For primes $p$, of the form $p\equiv 1 \mod 4$, also a different type of integral point sets occurs. To describe these sets, we need some new notation. For a prime $p\equiv 1\mod 4$, there is a unique element $\omega(p)\in\mathbb{N}$, with $\omega(p)<\frac{p}{2}$ and $\omega^2(p)\equiv -1\mod p$. By $\square_n=\{i^2\mid i\in\mathbb{Z}_n\}$ we denote the set of squares in $\mathbb{Z}_n$.

\begin{lemma}
  \label{lemma_ilig}
  For a prime $p\ge 3$, the set $\mathcal{P}=(1,\pm\omega(p))\cdot\square_p$ is a non-collinear integral point
  set over $\mathbb{Z}_p^2$ with cardinality $p$.
\end{lemma}
\begin{proof}
  For an odd prime $p$ we have exactly $\frac{p+1}{2}$ squares in $\mathbb{Z}_p$. Since $(0,0)$, $(1,\omega(p))$,
  and $(1,-\omega(p))$ are elements of $\mathcal{P}$, the point set is clearly non-collinear. For the property of
  pairwise integral distances we consider two arbitrary elements $q,q'\in\square_p$ and the corresponding distances
  \begin{eqnarray*}
    (q-q')^2+\omega^2(p)(q-q')^2&=& \overline{0},\\ 
    (q-q')^2+\omega^2(p)(q+q')^2&=& (2\omega(p))^2qq',\\
    (q+q')^2+\omega^2(p)(q-q')^2&=& 2^2qq',\\
    (q+q')^2+\omega^2(p)(q+q')^2&=& \overline{0}.\\
  \end{eqnarray*}
\end{proof}

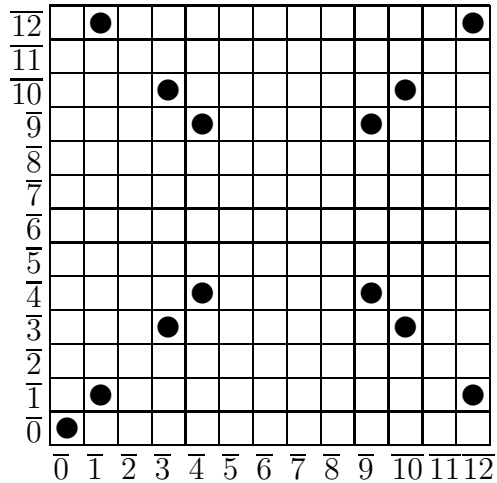
\begin{figure}[htp]
  \begin{center}
    \setlength{\unitlength}{0.45cm}
    \begin{picture}(14.2,14.2)
      \multiput(1.15,1)(0,1){14}{\line(1,0){13}}
      \multiput(1.15,1)(1,0){14}{\line(0,1){13}}
      \put(1.65,1.5){\circle*{0.6}}
      \put(2.65,2.5){\circle*{0.6}} 
      \put(4.65,4.5){\circle*{0.6}}
      \put(5.65,5.5){\circle*{0.6}}
      \put(10.65,10.5){\circle*{0.6}}
      \put(11.65,11.5){\circle*{0.6}}
      \put(13.65,13.5){\circle*{0.6}}
      \put(2.65,13.5){\circle*{0.6}}
      \put(4.65,11.5){\circle*{0.6}}
      \put(5.65,10.5){\circle*{0.6}}
      \put(10.65,5.5){\circle*{0.6}}
      \put(11.65,4.5){\circle*{0.6}}
      \put(13.65,2.5){\circle*{0.6}}
      \put(1.25,0){$\overline{0}$}
      \put(2.25,0){$\overline{1}$}
      \put(3.25,0){$\overline{2}$}
      \put(4.25,0){$\overline{3}$}
      \put(5.25,0){$\overline{4}$}
      \put(6.25,0){$\overline{5}$}
      \put(7.25,0){$\overline{6}$}
      \put(8.25,0){$\overline{7}$}
      \put(9.25,0){$\overline{8}$}
      \put(10.25,0){$\overline{9}$}
      \put(11.25,0){$\overline{10}$}
      \put(12.35,0){$\overline{11}$}
      \put(13.35,0){$\overline{12}$}
      \put(0.45,1.1){$\overline{0}$}
      \put(0.45,2.1){$\overline{1}$}
      \put(0.45,3.1){$\overline{2}$}
      \put(0.45,4.1){$\overline{3}$}
      \put(0.45,5.1){$\overline{4}$}
      \put(0.45,6.1){$\overline{5}$}
      \put(0.45,7.1){$\overline{6}$}
      \put(0.45,8.1){$\overline{7}$}
      \put(0.45,9.1){$\overline{8}$}
      \put(0.45,10.1){$\overline{9}$}
      \put(0,11.1){$\overline{10}$}
      \put(0,12.1){$\overline{11}$}
      \put(0,13.1){$\overline{12}$}
    \end{picture}
  \end{center}
  \caption{The integral point set $\Ps=(1,\omega(p))\cdot\square_p$ for $p=13$.}
  \label{fig_ilig}
\end{figure}

In Figure \ref{fig_ilig} we have depicted an integral point set, being constructed as described in Lemma \ref{lemma_ilig} for $p=13$. We remark that recently in \cite{finite_fields} it was proven, that integral point sets $\Ps$ over $\mathbb{Z}_p^2$, with cardinality $p\ge 3$, are either collinear or a translated version of the integral point set constructed in Lemma \ref{lemma_ilig}.

\subsection{Integral point sets over $\mathbf{\mathbb{Z}_n^2}$ with further conditions}

In the last subsection we have recognized, that integral point sets over $\mathbb{Z}_n^2$ are,
similar to integral point sets over $\mathbb{E}^2$, somewhat attracted by collinear sets. So we
investigate in this subsection integral point sets $\mathcal{P}$ over $\mathbb{Z}_n^2$, where no three points are collinear.

\begin{definition}
  By $\overline{\I}(n,\di)$ we denote the maximum number of points in semi-general position over $\Z$, where are
  pairwise distances are integral.
\end{definition}

\noindent
If we drop the condition of pairwise integral distances, our studied objects become very 
familiar discrete structures.
In the case of affine finite geometries  (classical \cite{Hi1} in the case of $\mathbb{Z}_n$ with 
$n$  a prime, Hjelmslev geometries \cite{0821.00012} in the other cases) point sets in semi-general position, 
with arbitrary pairwise distances, are called arcs in the case of planes or caps \cite{Br2} in the three
dimensional case. 
With the results from Subsection \ref{sub_sec_hom} in mind, we would like to mention the connection of 
these objects to linear coding theory, see i.e. \cite{Br1} for the details.

\noindent 
In Table \ref{table_n3ol} we give some values of $\overline{\mathcal{I}}(n,2)$ for small $n$, obtained by Algorithm \ref{alg_orderly} described later on.

\begin{table}[ht]
  \begin{center}
    \begin{tabular}{|r|r||r|r||r|r||r|r||r|r||r|r|}
      \hline $n$ & $\II$ & $n$ & $\II$ & $n$ & $\II$ & $n$ & $\II$ & $n$ & $\II$ & $n$ & $\II$ \\ 
      \hline   1 &     1 &  11 &     6 &  21 &     4 &  31 & 16 &  41 &     20 &  51 &      8 \\ 
      \hline   2 &     4 &  12 &     4 &  22 &     8 &  32 & 14 &  42 &      6 &  52 &     12 \\ 
      \hline   3 &     2 &  13 &     6 &  23 &    12 &  33 &  6 &  43 &     22 &  53 &     26 \\ 
      \hline   4 &     4 &  14 &     6 &  24 &     6 &  34 & 10 &  44 &     10 &  54 &$\ge 13$\\ 
      \hline   5 &     4 &  15 &     4 &  25 &    10 &  35 &  6 &  45 &     11 &  55 &      8 \\ 
      \hline   6 &     4 &  16 &     8 &  26 &    10 &  36 & 12 &  46 &     14 &  56 &     10 \\ 
      \hline   7 &     4 &  17 &     8 &  27 &    10 &  37 & 18 &  47 &     24 &  57 &     10 \\ 
      \hline   8 &     6 &  18 &    10 &  28 &     8 &  38 & 12 &  48 &      8 &  58 &$\ge 16$\\ 
      \hline   9 &     6 &  19 &    10 &  29 &    14 &  39 &  6 &  49 &$\ge 18$&  59 &     30 \\ 
      \hline  10 &     6 &  20 &     8 &  30 &     6 &  40 & 10 &  50 &$\ge 17$&  60 &      8 \\ 
      \hline 
    \end{tabular} 
    \caption{Values of $\overline{\mathcal{I}}(n,2)$ for small parameters $n$.}
    \label{table_n3ol}
  \end{center}
\end{table}

Now we want to derive an upper bound for $\overline{\I}(n,2)$, by relaxing the condition of pairwise integral distances. Let $\mathcal{P}$ be a point set over $\mathbb{Z}_n^2$ in semi-general position. We consider the lines $\left\{(i,j)\mid j\in\mathbb{Z}_n\right\}$ for $i\in\mathbb{Z}_n$. Since these $n$ lines form a partition of $\mathbb{Z}_n^2$ and each line can contain at most two points of $\mathcal{P}$, we obtain the trivial upper bound $\overline{\mathcal{I}}(n,2)\le 2n$. This is connected to a famous open problem in number theory \cite[sec. F4]{gu1}, where people work on an upper bound for the \textit{no-three-in-a-line} problem. Considering all lines  in $\mathbb{Z}_n^2$ we receive 
$$
  \overline{\mathcal{I}}(p,2)\le p+1
$$
for odd primes $p$ \cite{1005.51003} and 
$$
  \overline{\mathcal{I}}(n,2) \le n\cdot(1+ p^{- {\lceil {\frac{a+1}{2}} \rceil} } + p^{-a})
$$
where $p^a\mid n$ and $p^{a+1}\nmid n$ for a prime $p$ \cite{huizenga-max}.

Very recently for the case of odd primes $p$, tight bounds on $\overline{\I}(p,2)$ are proven \cite{finite_fields}:

\begin{theorem}
  \label{thm_semi_general}
  For $p\equiv 3\mod 4$ we have
  $$
    \overline{\I}(2,p)=\frac{p+1}{2}
  $$
  and for $p\equiv 1\mod 4$ we have
  $$
    \frac{p-1}{2}\le\overline{\I}(2,p)\le\frac{p+3}{2}.
  $$
\end{theorem}

We would like to remark that the known construction uses half of the points of the circle $\{(a,b)\in\mathbb{Z}_p^2\mid a^2+b^2=\overline{1}\}$, see \cite{finite_fields} for the details. For $p\equiv 1\mod 4$, $p\neq 5$ we conjecture $\overline{\I}(p,2)=\frac{p-1}{2}$.

\noindent
By a look at the situation in $\mathbb{E}^2$ and with the famous question of Erd\H{o}s in mind. it seems interesting to investigate integral point sets over $\mathbb{Z}_n^2$, where no three points are collinear and no four points are situated on a circle. 

\begin{definition}
  Four points $p_i=(x_i,y_i)$ in $\mathbb{Z}_n^2$ are said to be situated on a circle if there exist $a,b
  \in\mathbb{Z}_n$, $r \in\mathbb{Z}_n\backslash\{\overline{0}\}$ with
  $$
    (x_i-a)^2+(y_i-b)^2=r^2
  $$
  for all $i$.
\end{definition}

\noindent
We have the following necessary condition:

\begin{lemma}
  Four points $p_i=(x_i,y_i)$ in $\mathbb{Z}_n^2$ being situated on a circle fulfill
  \begin{equation}
    \label{eqn_circle}
    \left|\begin{array}{rrrr}
    x_1^2+y_1^2 & x_1 & y_1 & \overline{1} \\
    x_2^2+y_2^2 & x_2 & y_2 & \overline{1} \\
    x_3^2+y_3^2 & x_3 & y_3 & \overline{1} \\
    x_4^2+y_4^2 & x_4 & y_4 & \overline{1} \\
    \end{array}\right|=\overline{0}.
  \end{equation}
\end{lemma}

\begin{definition}
  By $\dot{\I}(n,\di)$ we denote the maximum number of points in $\Z$ with pairwise integral distances, where no
  three points are collinear and no four points are situated on a circle. Here we also talk of general position.
\end{definition}

\begin{table}[ht]
  \begin{center}
    \begin{tabular}{|r|r||r|r||r|r||r|r||r|r||r|r||r|r|}
      \hline  $n$ & $\IC$ & $n$ & $\IC$ & $n$ & $\IC$ & $n$ & $\IC$ & $n$ & $\IC$ & $n$ & $\IC$ & $n$ & $\IC$ \\ 
      \hline    1 &   1 &\!11\!&   4 &\!21\!&   4 &\!31\!&          6 &\!41\!&          9 &\!51\!&      7 &\!61\!& $\ge 9$\\ 
      \hline    2 &   4 &\!12\!&   4 &\!22\!&   8 &\!32\!&          8 &\!42\!&          6 &\!52\!& $\ge 9$&\!62\!&$\ge 11$\\ 
      \hline    3 &   2 &\!13\!&   5 &\!23\!&   5 &\!33\!&          4 &\!43\!&          8 &\!53\!& $\ge 9$&\!63\!&      8 \\ 
      \hline    4 &   4 &\!14\!&   6 &\!24\!&   4 &\!34\!&         10 &\!44\!&          8 &\!54\!&$\ge 11$&\!64\!&$\ge 10$\\ 
      \hline    5 &   4 &\!15\!&   4 &\!25\!&   6 &\!35\!&          5 &\!45\!&          8 &\!55\!&      6 &\!65\!&      7 \\ 
      \hline    6 &   4 &\!16\!&   6 &\!26\!&   8 &\!36\!&\!$\ge 10$\!&\!46\!&         10 &\!56\!&      6 &\!66\!&      8 \\ 
      \hline    7 &   3 &\!17\!&   5 &\!27\!&   7 &\!37\!&          7 &\!47\!&          7 &\!57\!&      6 &\!67\!& $\ge 9$\\ 
      \hline    8 &   4 &\!18\!&   8 &\!28\!&   6 &\!38\!&          8 &\!48\!&          8 &\!58\!&$\ge 11$&\!68\!&$\ge 10$\\ 
      \hline    9 &   4 &\!19\!&   5 &\!29\!&   7 &\!39\!&          6 &\!49\!&\!$\ge 11$\!&\!59\!& $\ge 9$&\!69\!&      7 \\ 
      \hline\!10\!&   6 &\!20\!&   6 &\!30\!&   6 &\!40\!&          6 &\!50\!&\!$\ge 12$\!&\!60\!&      8 &\!70\!& $\ge 9$\\ 
      \hline 
    \end{tabular} 
    \caption{Values of $\dot{\mathcal{I}}(n,2)$ for small parameters $n$.}
    \label{table_n4oc}
  \end{center}
\end{table}

\noindent
Trivially we have $\dot{\mathcal{I}}(n,2)\le\overline{\mathcal{I}}(n,2)$. 
In Table \ref{table_n4oc} we give some exact values of $\dot{\mathcal{I}}(n,2)$, obtained 
by Algorithm \ref{alg_orderly} described later on. One might conjecture that $\dot{\mathcal{I}}(n,2)$ is unbounded.

Because semi-general position or general position is a property of three or four points, respectively, we cannot apply our approach via clique search for the determination of $\overline{\I}(n,2)$ and $\dot{\I}(n,2)$ directly. Instead of going over to hypergraphs we use a variant of orderly generation \cite{winner}, which glues two integral point sets consisting of $r$ points, having $r-1$ points in common, to obtain recursively integral point sets of $r+1$ points. The used variant of orderly generation was introduced, and applied for the determination of the minimum distance $\dot{d}(n,2)$ of integral point sets in general position in $\mathbb{E}^2$, in \cite{phd_kurz,sascha-alfred}.

Now we go into detail. To describe integral point sets over $\mathbb{Z}_n^2$, we utilize the set $D_{n,2}$, 
where the coordinates of the points are \textit{reduced} with respect to the Lee weight via
$$
  \delta_n((x_1,y_1),(x_2,y_2))=\Big(\min(|\hat{x}_1-\hat{x}_2|,n-|\hat{x}_1-\hat{x}_2|),\min(|\hat{y}_1-
  \hat{y}_2|,n-|\hat{y}_1-\hat{y}_2|)\Big).
$$
By $\mathcal{B}=\{b_0,b_1,\dots,b_t\}$ we denote the subset of $D_{n,2}=\{\delta_n(0,x)\mid x\in\mathbb{Z}_n^2\}$, where the points $x$ are at integral distance to $0$. We define $b_0=(\overline{0},\overline{0})$. The numbering of the remaining $b_i$ is arbitrary but fix. Each integral point set $\Ps=\{p_1,\dots,p_r\}$ over $\mathbb{Z}_n^2$ is, up to translations and reflections, completely described by a matrix
$$
  \Delta_n(\Ps)=\Big(\iota\left(\delta_n\left(p_i,p_j\right)\right)\Big)_{i,j},
$$ 
where we set $\delta_n(p_i,p_i)=b_0$ and $\iota:\mathcal{B}\rightarrow\mathbb{N}$, $b_i\mapsto i$. We use these matrices as a data structure for integral point sets over $\mathbb{Z}_n^2$. Next we extend the natural order $\le$ on $\mathbb{N}$ to $\preceq$ for symmetric matrices, with zeros on the main diagonal as $\Delta_n$, by using a column-lexicographical order of the upper right matrix. A matrix $\Delta_n$ is said to be \emph{canonical} if $\Delta_n\ge \pi(\Delta_n)$ for every permutation $\pi\in S_r$ acting on the rows and columns of $\Delta_n$. If $\downarrow\!\!\!\Delta_n$ denotes the removal of the last column and last row of a matrix $\Delta_n$, then $\Delta_n$ is said to be \emph{semi-canonical} if $\downarrow\!\!\Delta_n\ge \downarrow\!\!\pi(\Delta_n)$ for 
every permutation $\pi\in S_r$. The function $\Gamma_r$ does the glueing of two integral point sets over 
$\mathbb{Z}_n^2$ consisting of $r$ points having $r-1$ points in common. The result of the function  $\Gamma_r$  
is an, with respect to $\preceq$,  ordered list of integral point sets consisting of $r+1$ points. By $\mathcal{L}_r$ we denote the ordered list of all semi-canonical matrices $\Delta_n$, with respect to $\preceq$, which correspond to integral point sets over $\mathbb{Z}_n^2$. It can be figured out easily that $\Gamma_r$ produces a list with at most two integral point sets. With these definitions we 
can state:    

\begin{algo}$\,$\\
  \label{alg_orderly}
  \textit{Input:} $\mathcal{L}_r$\\
  \textit{Output:} $\mathcal{L}_{r+1}$\\
  {\bf begin}\\
  \Tsm $\mathcal{L}_{r+1}=\emptyset$\\
  \Tsm{\bf loop over} $x_1\in\mathcal{L}_r$,\quad $x_1$ is canonical {\bf do}\\
  \Tsm\Tsm{\bf loop over} $x_2\in\mathcal{L}_r$,\quad $x_2\preceq x_1$, $\downarrow\!\!x_1=\downarrow\!\!x_2$ {\bf do}\\
  \Tsm\Tsm\Tsm{\bf loop over} $y\in\Gamma_r(x_1,x_2)$\\
  \Tsm\Tsm\Tsm\Tsm{\bf if} $y$ is semi-canonical {\bf then} add $y$ to $\mathcal{L}_{r+1}$ {\bf end}\\
  \Tsm\Tsm\Tsm{\bf end}\\
  \Tsm\Tsm{\bf end}\\
  \Tsm{\bf end}\\
  {\bf end}
\end{algo}

\noindent
A starting list $\mathcal{L}_3$ of the integral triangles can be generated by a nested loop. In order to apply Algorithm \ref{alg_orderly} for the determination of $\overline{\I}(n,2)$ or $\dot{\I}(n,2)$, we only have to modify it in that way, that it only accepts integral point sets in semi-general or general position, respectively, for the lists $\mathcal{L}_r$. 

\section{Integral point sets over $\mathbf{(\mathbb{R} / \mathbb{Z}n)^2}$}
In the previous section we have required also the coordinates of the point sets to be \textit{integral}.
This corresponds somewhat to integral point sets in $\mathbb{Z}^\di$. In this section we try to develop a 
setting for an analogous treatment of integral point sets in $\mathbb{E}^\di$ over the 
ring $\mathbb{Z}_n$ instead of $\mathbb{Z}$ for the distances. We start with $n=p$ being an odd prime.

Let $p$ be an odd prime, then $\mathbb{Z}_p$ is a finite field. Given three elements 
$a,b,c\in\mathbb{Z}_p\backslash\{\overline{0}\}$, which we consider as edge lengths of a triangle.
 Then we can determine a coordinate represention, given by three points
$(x_1,y_1)$, $(x_2,y_2)$, $(x_3,y_3)$ in $(\mathbb{R} / \mathbb{Z}p)^2$, as follows.
Due to translations, rotations and reflections we can assume $(x_1,y_1)=(\overline{0},\overline{0})$ and $(x_2,y_2)=(a,\overline{0})$. For the third point $(x_3,y_3)$ we get the system of equations
\begin{eqnarray*}
  x_3^2 + y_3^2 & = & b^2 , \\
  (x_3-a)^2 + y_3^2 & = & c^2 .
\end{eqnarray*}
Solving this system yields
\begin{eqnarray*}
  x_3 & = & \frac{b^2-c^2+a^2}{2a} ,\\
  y_3^2 & = & \frac{(a+b+c)(a+b-c)(a-b+c)(-a+b+c)}{(2a)^2} ,
\end{eqnarray*}
which is defined in $\mathbb{Z}_p$ because of $2a\neq\overline{0}$. By $\alpha(p)$ we denote the smallest quadratic non-residue in $\mathbb{Z}_p$. With the above system of equations it can be seen that $x_3\in\mathbb{Z}_p$ and $y_3$ is either also in $\mathbb{Z}_p$ or in $\mathbb{Z}_p\cdot \sqrt{\alpha(p)}$. Since this is similar to the case in $\mathbb{E}^\di$, see \cite{phd_kurz,paper_characteristic}, we define the characteristic of an integral triangle similarly.

\begin{definition}
  \label{def_characteristic}
  For an odd prime $p$ the characteristic of three side lengths $a,b,c\in\mathbb{Z}_p$ with
  $V^2=(a+b+c)(a+b-c)(a-b+c)(-a+b+c)\neq\overline{0}$ is defined as $\overline{1}$ if $V^2$ is a
  quadratic residue in $\mathbb{Z}_p$ and as $\alpha(p)$ otherwise.
\end{definition}

For the ease of notation we associate $\mathbb{E}_p^\di$ with $(\mathbb{R} / \mathbb{Z}p)^\di$.
We remark that the three points are collinear exactly if $V^2$ equals $\overline{0}$. So, similarly to the 
case in $\mathbb{E}^2$ \cite{54.0622.02}, we have the following lemma, where the determinant 
equals $V^2$, if we associate $a=\delta(v_1,v_2)$, $b=\delta(v_1,v_3)$, and $c=\delta(v_2,v_3)$.

\begin{lemma}
  Points $v_1,v_2,v_3\in\mathbb{E}_p^2$ are collinear if and only if their Euclidean distances 
  $\delta(v_i,v_j)$ fulfill
  $$
    \left|\begin{array}{rrrr}
      \delta^2(v_1,v_1) & \delta^2(v_1,v_2) & \delta^2(v_1,v_3) & \overline{1} \\
      \delta^2(v_2,v_1) & \delta^2(v_2,v_2) & \delta^2(v_2,v_3) & \overline{1} \\
      \delta^2(v_3,v_1) & \delta^2(v_3,v_2) & \delta^2(v_3,v_3) & \overline{1} \\
      \overline{1}      & \overline{1}      & \overline{1}      & \overline{0}
    \end{array}\right|=\overline{0}.
  $$
\end{lemma}

Our definition of the characteristic of an integral triangle in $\mathbb{Z}_p$ is properly chosen in the sense that we have the following theorem.

\begin{theorem}
  \label{thm_characteristic}
  In an integral point set over $\mathbb{E}_p^2$ where $p$ is an odd prime the characteristic of each
  non-degenerated triangle is equal.
\end{theorem}
\begin{proof}
  Without loss of generality we assume that the two triangles have two points in common and the points are given
  by the coordinates $(\overline{0},\overline{0})$, $(\overline{0},a)$, $(x,y\sqrt{c})$, $(x',y'\sqrt{c'})$, where
  $a$, $x$, $x'$, $y$, $y'$ are elements of $\mathbb{Z}_p$ and $c,c'$ are the characteristics. The squared distance
  of the last two points is given by
  $$
    (x-x')^2+(y\sqrt{c}-y'\sqrt{c'})^2=(x-x')^2+y^2c-2yy'\sqrt{cc'}+{y'}^2c'.
  $$
  Because this number must be an element of $\mathbb{Z}_p$ we have that $cc'$ is a quadratic residue in
  $\mathbb{Z}_p$ yielding $c=c'$.
\end{proof}

As we have proceeded completely analogous to the case in $\mathbb{E}^\di$ we can generalize Definition \ref{def_characteristic} and Theorem \ref{thm_characteristic}.

\begin{definition}
  \label{def_characteristic_general}
  For an odd prime $p$ the characteristic of an integral point set with $\di+1$ points in $\mathbb{E}_p^\di$ given by 
  its distances $\delta_{i,j}$ is $1$ if $V^2_\di$ is a quadratic residue in $\mathbb{Z}_p$ and $\alpha(p)$ otherwise, where 
  $$
    V^2_\di=\left|\begin{array}{rrrr}
      \delta^2_{1,1} & \dots & \delta^2_{1,\di+1} & \overline{1}\\
      \vdots & \ddots & \ddots & \vdots \\
      \delta^2_{\di+1,1} & \dots & \delta^2_{\di+1,\di+1} & \overline{1}\\
      \overline{1} & \dots & \overline{1} & \overline{0}  
    \end{array}\right|.
  $$
\end{definition}

\begin{theorem}
  In an integral pointset over $\mathbb{E}_p^\di$ where $p$ is an odd prime the characteristic of each non-degenerated
  simplex the same.
\end{theorem}
\begin{proof}
  We do the corresponding calculations as in \cite{paper_characteristic} over $\mathbb{Z}_p$ instead of $\mathbb{Q}$.
\end{proof}

For completeness we give a necessary coordinatefree criterion for $\di+2$ points being situated on an $\di$-dimensional sphere.

\begin{lemma}
  \label{lemma_sphere_general}
  If $\di+2$ points in $\mathbb{E}_n^\di$ described by their distances $\delta_{i,j}$ are situated on 
  an $\di$-dimensional sphere then 
  $$
    \left|\begin{array}{rrr}
      \delta^2_{1,1} & \dots & \delta^2_{1,\di+1} \\
      \vdots & \ddots & \vdots \\
      \delta^2_{\di+1,1} & \dots & \delta^2_{\di+1,\di+1}
    \end{array}\right|= \overline{0}.
  $$
\end{lemma}

So far we have transferred the theory of integral point sets in $\mathbb{E}^\di$ to integral point sets over $\mathbb{E}_p^\di$ for odd primes $p$. For general $n$ instead of $p$ there are some twists if we use coordinates. The most natural approach to settle these would be, with respect to the situation in $\mathbb{E}^\di$, to leave out coordinates and use 
Mengers characterization of embedable distance matrices \cite{54.0622.02} and replace the conditions over $\mathbb{Z}$ by conditions over $\mathbb{Z}_n$.

\begin{definition}
  \label{def_integral_pointset}
  An integral point set $\Ps$ over $\mathbb{E}_n^\di$ is a set of $r\ge\di+1$ points with distances 
  $\delta_{i,j}\in\mathbb{Z}_n\backslash\{\overline{0}\}$ for $1\le i\neq j\le r$ which fulfill 
  $$
    V_{t-1}^2(\{i_1,\dots,i_t\})=\left|\begin{array}{rrrr}
      \delta^2_{i_1,i_1} & \dots & \delta^2_{i_1,i_t} & \overline{1}\\
      \vdots & \ddots & \ddots & \vdots \\
      \delta^2_{i_t,i_1} & \dots & \delta^2_{i_t,i_t} & \overline{1}\\
      \overline{1} & \dots & \overline{1} & \overline{0}  
    \end{array}\right|=\overline{0}
  $$
  for each subset of points $\{i_1,\dots,i_t\}$ of cardinality $t=\di+2$ and $t=\di+3$, and there exists a subset 
  $\{\tilde{\text{\i}}_1,\dots,\tilde{\text{\i}}_t\}$ of cardinality $t=\di+1$ with 
  $V_{t-1}^2(\{\tilde{\text{\i}}_1,\dots,\tilde{\text{\i}}_t\})\neq \overline{0}$.
\end{definition}

To model the extra conditions we could define that $\Ps$ is in semi-general position if for every $\di+1$ points $\{i_1,\dots,i_{\di+1}\}$ we have $V_{\di+1}^2(\{i_1,\dots,i_{\di+1}\})\neq\overline{0}$ and that $\Ps$ is in general position if the condition of Lemma \ref{lemma_sphere_general} is fulfilled. We remark that for $\di=2$ the determinant of Lemma \ref{lemma_sphere_general} can be factorized to
\begin{eqnarray*}
  -(\delta_{1,2}\delta_{3,4}+\delta_{1,3}\delta_{2,4}+\delta_{1,4}\delta_{2,3})
  (\delta_{1,2}\delta_{3,4}+\delta_{1,3}\delta_{2,4}-\delta_{1,4}\delta_{2,3})\cdot\\
  (\delta_{1,2}\delta_{3,4}-\delta_{1,3}\delta_{2,4}+\delta_{1,4}\delta_{2,3})
  (-\delta_{1,2}\delta_{3,4}+\delta_{1,3}\delta_{2,4}+\delta_{1,4}\delta_{2,3}).
\end{eqnarray*}
For $\di=2$ we also have
\begin{eqnarray*}
  V_2^2(\{1,2,3\})&=&(\delta_{1,2}+\delta_{1,3}+\delta_{2,3})(\delta_{1,2}+\delta_{1,3}-\delta_{2,3})\cdot\\
  &&(\delta_{1,2}-\delta_{1,3}+\delta_{2,3})(-\delta_{1,2}+\delta_{1,3}+\delta_{2,3}).
\end{eqnarray*}
So one may leave out the first factor and request that one of the remaining factors equals $\overline{0}$ instead of the condition in Definition \ref{def_integral_pointset} and the condition in Lemma \ref{lemma_sphere_general}, respectively. For $\di\ge 3$ the two corresponding determinants are irreducible \cite{irreducible}.


\noindent
Another way to generalize integral point sets is to consider the edge lengths and coordinates as elements in a finite field $\mathbb{F}_{p^k}$ or a commutative ring $\mathcal{R}$ instead of $\mathbb{F}_p=\mathbb{Z}_p$. For some results we refer to \cite{algo,finite_fields}. Here we only give a very general definition of an integral point set over an commutative ring $\mathcal{R}$:

\begin{definition}
  For a commutative ring $\mathcal{R}$ a set $\mathcal{P}$ of $n$ points in $\mathcal{R}^\di$ is called an integral 
  point set if for each $(x_1,\dots,x_\di),(y_1,\dots,y_\di)\,\in\mathcal{R}^\di$ there exists an element
  $d\in\mathcal{R}$ fulfilling
  $$
    \sum_{i=1}^\di (x_i-y_i)^2\,=\,d^2.
  $$
\end{definition}

\section{Conclusion}
We have generalized the theory of integral point sets over $\mathbb{Z}^m$ to integral point sets over $\mathbb{Z}_n^\di$. Some exact values $\mathcal{I}(n,\di)$ of the maximal cardinality of a set with pairwise integral distances in $\mathbb{Z}_n^\di$ with or without further conditions on the position are given together with algorithms to determine them.

There are two connections to coding theory, first via the special case of arcs and caps, secondly by the observation  that
$\mathcal{I}(n,\di)$ leads to a class of codes where the Hamming distances of the codewords have to fulfill certain modular restrictions. 

For odd primes $p$ the theory of integral point sets in $\mathbb{E}^\di$ is transferred to a theory of integral point sets over $\mathbb{E}_p^\di$ including the fundamental theorem about the characteristic of an integral simplex. 

There are some open questions left and the given results motivate for further research on integral point sets over $\mathbb{Z}_n^\di$ and $\mathbb{E}_n^\di$, as they seem to be interesting combinatorial structures.

\bibliography{ipozn}
\bibdata{ipozn}
\bibliographystyle{plain}

\end{document}